\documentclass[reqno]{amsart}
\usepackage{graphicx}
\usepackage{epsfig}
\usepackage{subfigure}
\theoremstyle{proposition}

\theoremstyle{definition}

\theoremstyle{remark}

\newcommand{\g}{\mathcal{G}}
\newcommand{\f}{\mathcal{F}}
\newcommand{\me}{\mathcal{E}}

\newcommand{\beq}{\begin{equation}}
\newcommand{\eeq}{\end{equation}}

\newcommand{\be}{\begin{enumerate}}
\newcommand{\ee}{\end{enumerate}}
\newcommand{\bi}{\begin{itemize}}
\newcommand{\ei}{\end{itemize}}
\newcommand{\bd}{\begin{description}}
\newcommand{\ed}{\end{description}}

\begin{document}

\title{Measuring the "non-stopping timeness" of ends of previsible sets}


\author{Ju-Yi Yen$^{(1),(2)}$}
\address{(1) Vanderbilt University, Nashville, Tennessee 37240, USA}
\email{ju-yi.yen@vanderbilt.edu}
\address{ (2) Institute of Mathematics, Academia Sinica, Taipei, Taiwan}
\curraddr{}
\thanks{}
\author{Marc Yor$^{(3),(4)}$}
\address{(3) Laboratoire de Probabilit\'{e}s et Mod\`{e}les Al\'{e}atoires,
Universit\'{e} Pierre et Marie Curie, Case Courrier 188, 4, Place
Jussieu, 75252 Paris, Cedex 05,  France}
\address{(4) Institut Universitaire de France} \email{}
\thanks{}


\date{October 3, 2008}

\dedicatory{}

\begin{abstract}
In this paper, we propose several "measurements" of the
"non-stopping timeness" of ends $\g$ of previsible sets, such that
$\g$ avoids stopping times, in an ambiant filtration. We then study
several explicit examples, involving last passage times of some
remarkable martingales.
\end{abstract}

\maketitle

\section{Introduction: About ends of previsible sets}
In this paper, we are interested in random times $\g$ defined on a
filtered probability space $(\Omega,\mathcal{F},(\mathcal{F}_t),P)$
as ends of $(\mathcal{F}_t)$ previsible sets $\Gamma$, that is:
\begin{equation}
\g\equiv \g_\Gamma = \sup\{t: \ (t,\omega) \in \Gamma\}
\end{equation}
For simplicity, we shall make the following assumptions:
\begin{description}
\item[(C)] All $((\f_t),P)$ martingales are continuous; \item[(A)]
For any $(\f_t)$ stopping time $T$, $P(\g=T)=0$.
\end{description}
[(C) stands for "continuous", and (A) for "avoiding"]. \\
To such a random time, one associates the Az\'{e}ma supermartingale:

\beq Z_t^\g=P(\g>t|\f_t), \eeq which, under (CA), admits a
continuous version.

In a number of questions, it is very interesting to consider the
smallest filtration $(\f_t^\prime)_{t\ge 0}$, which contains
$(\f_t)$, and makes $\g$ a stopping time; this filtration is usually
denoted $(\f_t^\g)_{t\ge 0}$. One of the interests of $(Z_t^\g)$ is
that it allows to write any $(\f_t)$ martingale as a semimartingale
in $(\f_t^\g)_{t\ge 0}$; see e.g. \cite{Jeulin, Jeulinyor,
Mansuyyor, N}, for both general formulae and many examples.

Recently, it has been understood that Black-Scholes like formulae
are closely related with certain such $\g$'s, thus throwing a new
light on a cornerstone of Mathematical Finance, see, e.g. \cite{MRY,
MRY1}. In the present paper, with (A) as our essential hypothesis,
we would like to measure "how much $\g$ differs from a $(\f_t)$
stopping time". The remainder of this paper consists in two
sections: - In Section 2, we propose several criterions to measure
the NST ($\equiv$ Non Stopping Timeness) of $\g$'s which satisfy
(CA); it is not surprising that we get interested in the function:
$$ m_\g(t)=E\left[\left(1_{(\g\ge t)}-P(\g\ge
t|\f_t)\right)^2\right].$$ - In Section 3, we compute explicitly
this function $m_\g$ for various examples, where $\g$ is the last
passage time at a level of a martingale which converges to 0, as
$t\rightarrow \infty$.

\section{Several possible "NST" criterions}
\subsection{A fundamental function: $m_\g(t)$}

Let us consider several possible criterions: first of all, if $\g$
were a $(\f_t)$ stopping time, then the process $1_{(\g\ge t)}$
would be identically equal to $Z_t^\g \equiv P(\g\ge t|\f_t)$. Thus,
the function $(m_\g(t), \ t \ge 0)$ tells us about the NST of $\g$.
Next, a simple but useful remark is that:
\begin{equation}
m_\g(t) = E\left[Z_t^\g\left(1-Z_t^\g\right)\right].
\end{equation}

\subsection{Definition of $m_\g^*$, $m_\g^{**}$, and $\widetilde{m}_\g$ and
proof that $m_\g^{**}=1/4 = \widetilde{m}_\g $} Instead of
considering the "full"function $(m_\g(t), \ t\ge 0)$, we may
consider only: 
\beq m_\g^*=\underset{t\ge 0}{\sup} \ m_\g(t) \eeq
as a "global" measurement of the NST of $\g$.

Here are two other, a priori natural, measurements of the NST of
$\g$:
\beq m_\g^{**} = E\left[\underset{t\ge 0}{\sup}
\left(Z_t^\g\left(1-Z_t^\g\right)\right)\right] \eeq and
 \beq \widetilde{m}_\g =
\underset{T\ge 0}{\sup} \ E\left[Z_T^\g\left(1-Z_T^\g\right)\right]
\eeq
where $T$ runs over all $(\f_t)$ stopping times.

However, we cannot expect to learn very much from $m_\g^{**}$ and
$\widetilde{m}_\g$, since it is easily shown that
\beq
 m_\g^{**}=1/4 = \widetilde{m}_\g.
\eeq
\begin{proof}
\mbox{} \bi
\item[(i)]
The fact that $m^{**}_\g=1/4$ follows immediately from: \\
$\underset{x\in[0,1]}{\sup} \ (x(1-x))=1/4$, and the fact that,
a.s., the range of the process $(Z_t^\g, t\ge 0)$ is $[0,1]$ since
$Z_0^\g=1$, $Z_\infty^\g=0$, and $(Z_t^\g, t\ge 0)$ is continuous.
\item[(ii)] Let us consider $T_a=\inf \{t: Z_t^\g=a\}$, for $0<a<1$. Then: \\
$Z_t^\g(1-Z_t^\g)\mid_{t=T_a}=a(1-a)$; hence, \\
$\underset{a\in]0,1[}{\sup}\left[Z_{T_a}^\g\left(1-Z_{T_a}^\g\right)\right]
=\underset{a\in]0,1[}{\sup}\left(a(1-a)\right)=1/4.$ \ei
\end{proof}

\subsection{The optional stopping time discrepancy $\mu_\g$}
Let us not forget either the nice characterization \cite{KM}, of
stopping times, among random times, as the times $\tau$ such that
for every bounded martingale $(M_t)_{t \ge 0}$ one has $$M_\tau=
E\left[M_\infty|\f_\tau\right]$$ where, here under our hypothesis
(C), we may define $\mathcal{F}_\tau=\sigma\{H_\tau; \
H \ \mbox{previsible}\}$. Thus, as a 4th possible 
measurement of the NST of $\g$, we may take: \beq \mu_\g =
\underset{\underset{E(M_\infty^2)\le 1}{M_\infty \in
L^2(\f_\infty)}}{\sup}E\left[\left(M_\g-E\left[M_\infty|\f_\g\right]\right)^2\right].\eeq

\subsection{Distance from stopping times}
We introduce:

$$\nu_\g=\underset{T\ge 0}{\inf}E\left[|\g-T|\right]$$ where $T$ runs over all $(\f_t)$ stopping times.
However, this quantity may be infinite as $\g$ may have infinite
expectation. A more adequate distance may be:

$$\nu_\g^\prime=\underset{T\ge 0}{\inf}\left(E\left[\frac{|\g-T|}{1+|\g-T|}\right]\right)$$
We note that this distance was precisely computed by du
Toit-Peskir-Shiryaev in the example they consider \cite{dTPS}.

It would be nice to be able to estimate $\mu_\g$ and/or $\nu_\g$,
$\nu_\g^\prime$ in terms of the function $m_\g(t)$, but we have not
been able to obtain any result in this direction. We shall now
concentrate uniquely on the study of $(m_\g(t), \ t\ge 0)$.
\section{A study of several interesting examples of functions $m_\g(t)$}
\subsection{Some general formulae}
We shall compute $(m_\g(t), \ t\ge 0)$ in some particular cases
where:
$$\g=\g_K=\sup\{t\ge 0: M_t=K\}, \ K\le 1,$$ with
$M_0=1$, $M_t\ge 0$, a continuous local martingale such that
$M_t\underset{t\rightarrow \infty}{\longrightarrow}0$, we recall
that (see, e.g. \cite{Jeulin,Mansuyyor}):
$$Z_t=P(\g_K\ge t|\f_t)=1\wedge \left(\frac{M_t}{K}\right)$$ thus
\beq \label{mt}
m(t)=E\left[Z_t\left(1-Z_t\right)\right]=\frac{1}{K^2}E\left[
M_t\left(K-M_t\right)^+\right] \eeq

\subsection{The particular case $M_t=\me_t=\exp(B_t-t/2)$,
with $(B_t)$ a standard Brownian motion, and
$\g_K=\sup\{t:\me_t=K\}$, $(K\le1)$}
From formula~(\ref{mt}), we deduce:
\begin{eqnarray*}
m_K(t)&=& \frac{1}{K^2}E\left[\me_t\left(K-\me_t\right)^+\right]\\
&=&
\frac{1}{K^2}E\left[\left(K-\exp\left(B_t+\frac{t}{2}\right)\right)^+\right]
\ \
\mbox{(by Cameron-Martin)} \nonumber \\
&=& \frac{1}{K^2}\left\{KP\left(\exp\left(B_t+\frac{t}{2}\right)<K\right) \right.\\
&& \ \ \ \left.
-E\left[1_{\left(\exp\left(B_t+\frac{t}{2}\right)<K\right)}\exp\left(B_t+\frac{t}{2}\right)\right]
\right\}\\
(K=e^l)&=& e^{-l}P\left(B_t+\frac{t}{2}<l\right)-e^te^{-2l}P\left(B_t+\frac{3t}{2}<l\right)\\
&=& P\left(B_1<-\frac{3\sqrt{t}}{2}+\frac{l}{\sqrt{t}}\right)\left(e^{-l}-e^{t-2l}\right)\\
&& \ \ \ \ +
e^{-l}P\left(-\frac{3\sqrt{t}}{2}+\frac{l}{\sqrt{t}}<B_1<-\frac{\sqrt{t}}{2}+\frac{l}{\sqrt{t}}\right).
\end{eqnarray*}

Particularizing again, for $M_t=\me_t=\exp(B_t-t/2)$, and $K=1$,
then
\begin{eqnarray}
m(t)&=& E\left[\me_t\left(1-\me_t\right)^+\right] \nonumber \\
&=&
P\left(B_1<-\frac{3\sqrt{t}}{2}\right)\left(1-e^t\right)+P\left(-\frac{3\sqrt{t}}{2}
<B_1<-\frac{\sqrt{t}}{2}\right). \label{mt1}
\end{eqnarray}
Figure 1 presents the graphs of $m_K(t)$ for some $K$'s.

\begin{figure}[htp]\label{mtkplot}
\includegraphics[width=4in,height=6.4in]{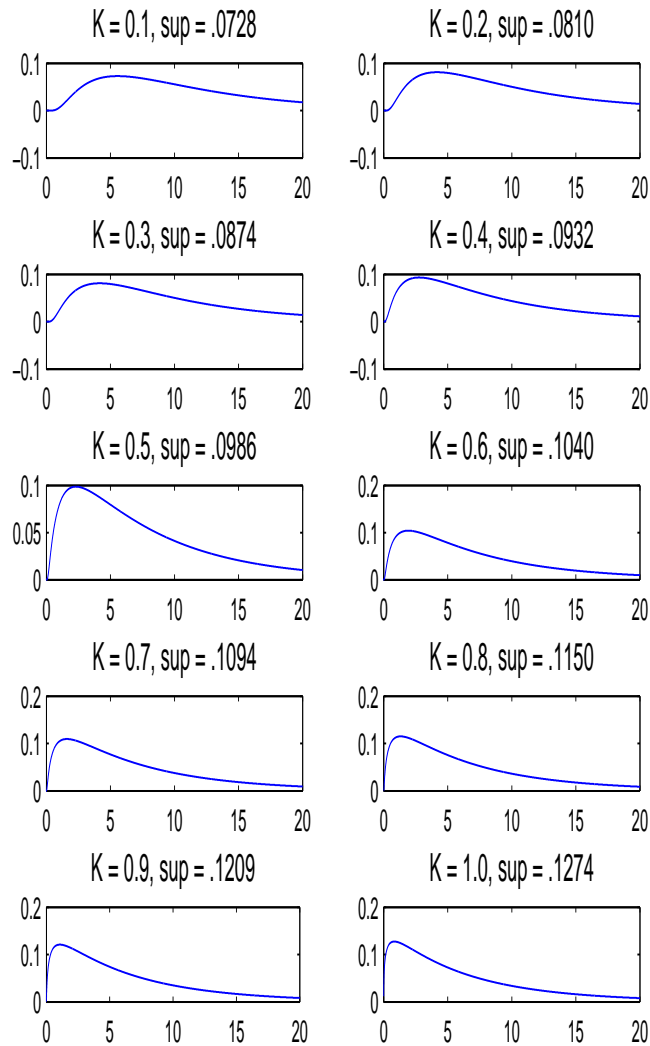}
\caption{\footnotesize\sf Graphs of $m_K(t)$, for $K=0.1,0.2,...1.$}
\end{figure}

\noindent\underline{Comments}: It seems that
$m_K^*=\underset{t\ge0}{\sup} \ (m_K(t))$ increases with $K$, at
least for $K\le 1$. Can this be proven? If so, it indicates that, as
$K$ varies, the $\g_K$'s are increasingly further from stopping
times.


\subsection{The case $\g=\g_{T_a}=\sup\{t<T_a:B_t=0\}$}

From line 4 of Table $(1\alpha)$ of Progressive Enlargements, p.32
of \cite{Mansuyyor}, we obtain:
$$Z_t=1-\frac{1}{a}B^+_{t\wedge T_a}$$ Thus, we obtain:

\begin{eqnarray}
m_\g(t)&=&E\left[\left(\frac{1}{a}B_{t\wedge T_a}^+\right)\left(1-\frac{1}{a}B_{t\wedge T_a}^+\right)\right] \nonumber \\
&=&\frac{1}{a^2}E\left[1_{(t<T_a)}1_{(B_t>0)}B_t(a-B_t)\right] \nonumber \\
&=& \frac{1}{a^2}E\left[1_{(S_t<a)}1_{(B_t>0)}B_t(a-B_t)\right] \nonumber \\
&=&
\frac{1}{a^2}E\left[\left(S_1<\frac{a}{\sqrt{t}}\right)1_{(B_1>0)}tB_1\left(\frac{a}{\sqrt{t}}-B_1\right)
\right] \nonumber \\
&=& \frac{1}{x^2}\varphi(x), \ \ \mbox{where:   $x=\frac{a}{\sqrt{t}}$   and:} \nonumber \\
\varphi(x)&=& E\left[1_{(S_1<x)}1_{(B_1>0)}B_1(x-B_1)\right]
\end{eqnarray}
Now, it remains to compute the function $\varphi$. We note that
$$\varphi(x)=E\left[B_1^+(x-B_1)^+\right]-E\left[1_{(S_1>x)}B_1^+(x-B_1)^+\right].$$
We shall take advantage of the very useful formula:
$$P(S_1>x|B_1=a)=\exp(-2x(x-a)), \ x\ge a>0.$$ (This formula is easily seen to be equivalent to the well-known expression of the joint density of $(S_1,B_1)$; see, e.g., \cite{KS}, p.425.)
Thus, we find:
\begin{eqnarray*}
\varphi(x) &=& \frac{1}{\sqrt{2\pi}}\int_{0}^{x}dy \ y(x-y)\left(\exp\left(-\frac{y^2}{2}\right)-\exp\left(-\frac{1}{2}(2x-y)^2\right)\right)\\
&\equiv& \frac{x^3}{\sqrt{2\pi}}\int_{0}^{1}du \
u(1-u)\left(\exp\left(-\frac{x^2u^2}{2}\right)-\exp\left(-\frac{x^2}{2}(2-u)^2\right)\right)
\end{eqnarray*}
Thus:
$$\frac{\varphi(x)}{x^2}=\frac{x}{\sqrt{2\pi}}\int_{0}^{1}du \ u(1-u)\left(\exp\left(-\frac{x^2u^2}{2}\right)-\exp\left(-\frac{x^2}{2}(2-u)^2\right)\right).$$

\subsection{The case $\g=L_a=\sup\{u:R_u=a\}$}

From line 6 of Table $(1\alpha)$ of Progressive Enlargements, p.32
of \cite{Mansuyyor}, we obtain:
$$Z_t=1\wedge\left(\frac{a}{R_t}\right)^{2\mu}.$$
Here, $(R_u)$ is $BES_0(d), \ d=2(\mu+1)$. For that second example,
we get:
\begin{eqnarray*}
m_\g(t)&=&
E\left[\left(1\wedge\left(\frac{a}{R_t}\right)^{2\mu}\right)
\left(1-1\wedge\left(\frac{a}{R_t}\right)^{2\mu}\right)\right]\\
&=&
E\left[1\left(\frac{a}{\sqrt{t}R_1}<1\right)\left(\frac{a}{\sqrt{t}R_1}\right)^{2\mu}
\left(1-\left(\frac{a}{\sqrt{t}R_1}\right)^{2\mu}\right)\right]\\
&=& \varphi_\mu\left(\frac{a^2}{2t}\right)
\end{eqnarray*}
where using the fact that
$R^2_1\stackrel{\rm{(law)}}{=}2\gamma_{d/2}$, we get (recall:
$\frac{d}{2}=\mu+1$):
$$\varphi_\mu(z)=\frac{1}{\Gamma(\mu+1)}\left\{z^\mu e^{-z}-z^{2\mu}\int_{z}^{\infty}
\frac{du}{u^\mu}e^{-u}\right\}.$$

\begin{figure}[htp]\label{mtkplot}
\includegraphics[width=3.2in]{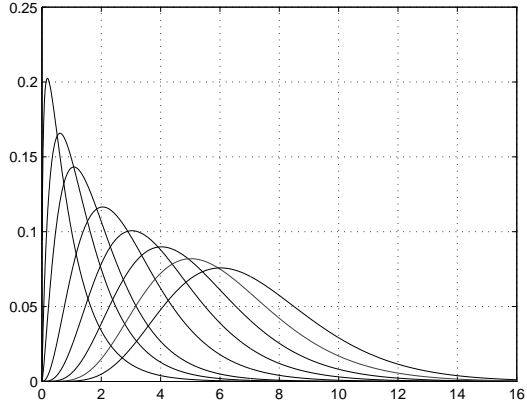}
\caption{\footnotesize\sf Graphs of $\varphi_\mu(z)$, for
$\mu=1/2,1,3/2,5/2,7/2,9/2,11/2,13/2$, and that
$z_{1/2}=0.19, z_{1}=0.61, z_{3/2}=1.08, z_{5/2}=2.05, z_{7/2}=3.04,
z_{9/2}=4.03, z_{11/2}=5.02, z_{13/2}=6.02$.}
\end{figure}

Figure 2 presents the graphs of $\varphi_\mu$ for
$\mu=1/2,1,3/2,5/2,7/2,9/2,11/2$ and $13/2$. We also approximate
$z_\mu$, the unique $>0$ real which achieves the max of
$\varphi_\mu$. This will give us the value
$m_\mu\stackrel{\rm{def}}{=}m_\g^*$, for these $\g\equiv
\mathcal{L}_a$ (note that, for a given $\mu$, the value does not
depend on $a$; this is because of the scaling property).

It is not difficult to show that: $z_\mu$ is the unique solution of
$$(E_\mu): \ \ \frac{1}{2z}=
\int_{0}^{\infty}\frac{dh}{(1+h)^\mu}e^{-hz}$$ and also
$$m_\mu=\frac{1}{\Gamma(\mu+1)}e^{-z_\mu}\frac{(z_\mu)^\mu}{2}.$$
Note that
$$m_\mu \le m_\mu^\prime\stackrel{\rm{def}}{=}
\frac{1}{\Gamma(\mu+1)}\ \underset{z\ge 0}{\sup} \
\left(e^{-z}\frac{z^\mu}{2}\right).$$
Figure 3 presents the graphs of $m_\mu$ and $m_\mu^\prime$.




\begin{figure}[ht]
\centering \subfigure{
\includegraphics[width=2.38in]{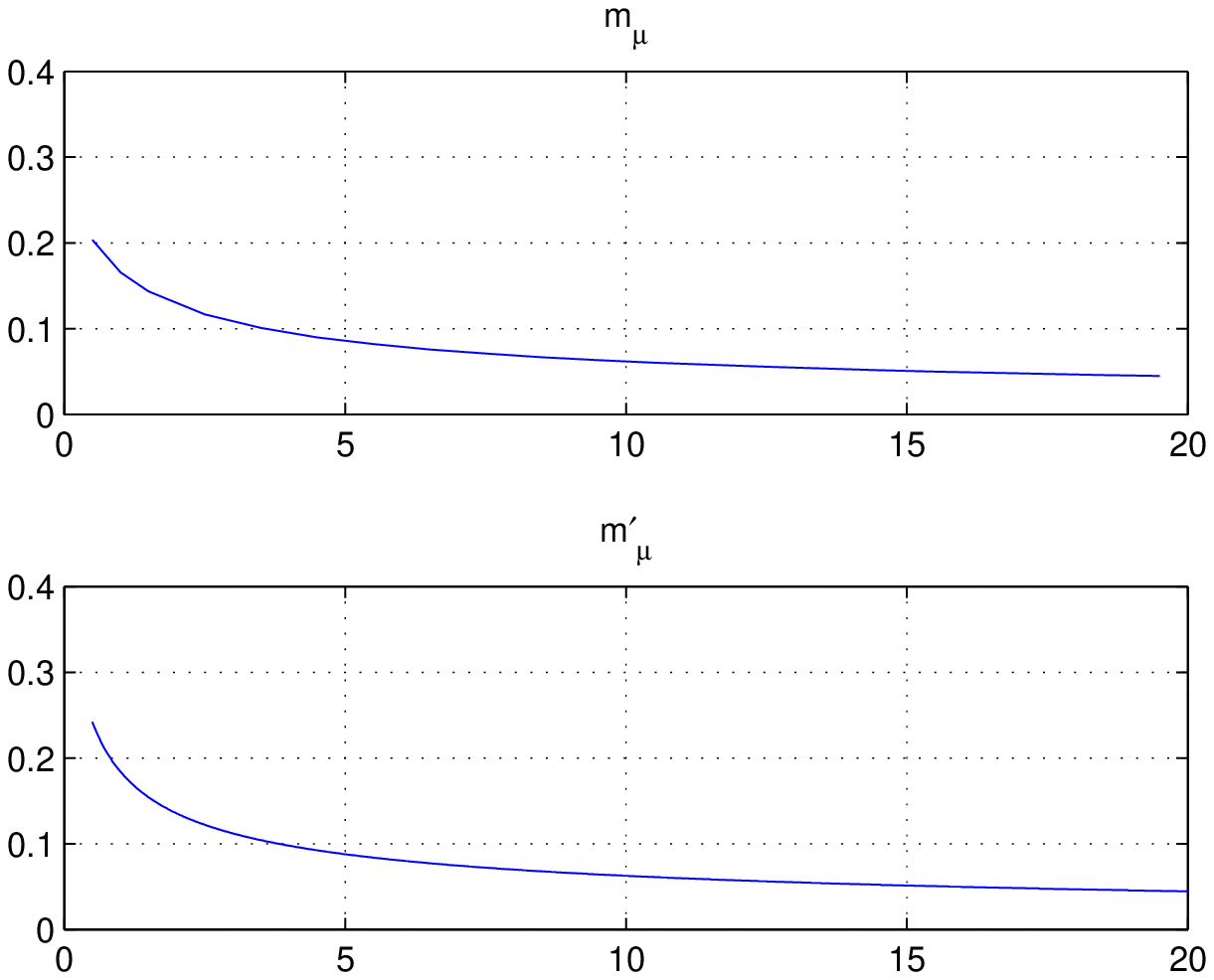}
\label{fig:subfig1} } \subfigure{
\includegraphics[width=2.35in]{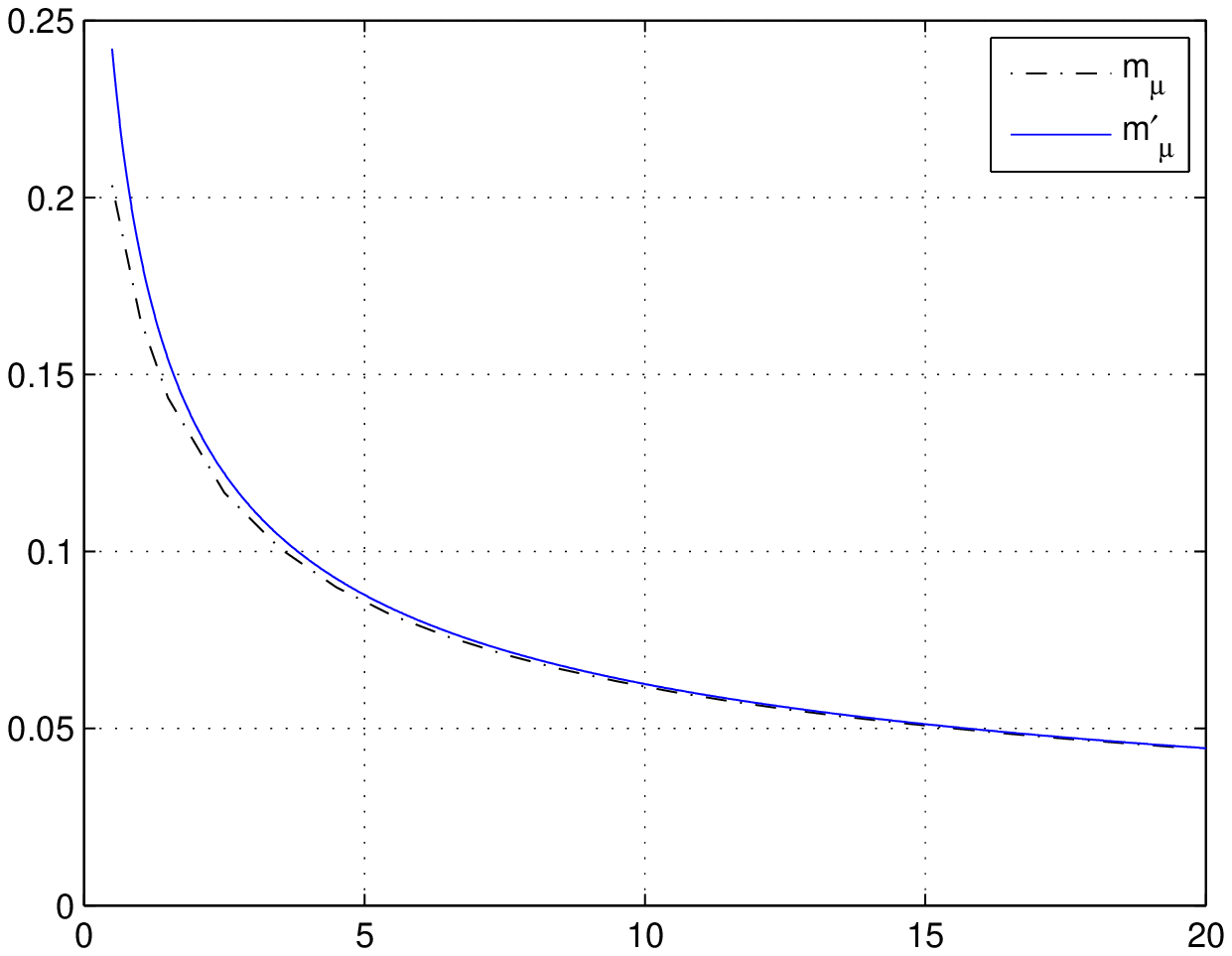}
\label{fig:subfig2} }  \caption{Graphs of $m_\mu$ and
$m_\mu^\prime$}
\end{figure}



\end{document}